\documentclass[12pt]{amsart}

\usepackage{geometry}                
\usepackage{graphicx,graphics}
\usepackage{amssymb,enumerate}
\usepackage{amsmath,amsfonts,amsthm}\usepackage{exscale}
\usepackage{latexsym}\usepackage{pict2e}
\usepackage{verbatim}\usepackage[all,poly,curve,knot,arrow]{xy}
\usepackage{fancybox}
\usepackage{wasysym}

\theoremstyle{plain}
\newtheorem{thm}{Theorem}[section]

\newtheorem{lem}[thm]{Lemma}
\newtheorem{prop}[thm]{Proposition}
\newtheorem{cor}[thm]{Corollary}

\newtheorem{ques}[thm]{Question}

\theoremstyle{definition}
\newtheorem{ob}[thm]{Observation}

\newtheorem{ex}[thm]{Example}
\newtheorem{defn}[thm]{Definition}
\newtheorem{rmk}[thm]{Remark}

\newtheorem{fact}[thm]{Fact}

\newcommand{\fr}{\frac}
\newcommand{\we}{\wedge}
\newcommand{\lf}{\lfloor}
\newcommand{\rf}{\rfloor}

\newcommand{\be}{\begin{enumerate}}
\newcommand{\ee}{\end{enumerate}}
\newcommand{\bd}{\begin{defn}}
\newcommand{\ed}{\end{defn}}
\newcommand{\bp}{\begin{prop}}
\newcommand{\ep}{\end{prop}}
\newcommand{\bl}{\begin{lem}}
\newcommand{\el}{\end{lem}}
\newcommand{\bq}{\begin{ques}}
\newcommand{\eq}{\end{ques}}

\newcommand{\bc}{\begin{cor}}
\newcommand{\ec}{\end{cor}}

\newcommand{\x}{\texttt{x}}
\newcommand{\lam}{\lambda}

\newcommand{\bt}{\begin{thm}}
\newcommand{\et}{\end{thm}}
\newcommand{\bpf}{\begin{proof}}
\newcommand{\epf}{\end{proof}}
\newcommand{\bex}{\begin{ex}}
\newcommand{\eex}{\end{ex}}
\newcommand{\bft}{\begin{fact}}
\newcommand{\eft}{\end{fact}}
\newcommand{\brk}{\begin{rmk}}
\newcommand{\erk}{\end{rmk}}
\newcommand{\al}{\alpha}\newcommand{\gam}{\gamma}

\newcommand{\rr}{\mathbb{R}}

\newcommand{\ba}{\begin{align*}}
\newcommand{\ea}{\end{align*}}
\newcommand{\tn}{\textnormal}

\newcommand{\bit}{\begin{itemize}}
\newcommand{\eit}{\end{itemize}}

\newcommand{\bcm}{}
\newcommand{\cref}[1]{(\ref{#1})}

\newcommand{\ci}{\CIRCLE}
\newcommand{\cci}{\Circle}

\newcommand{\di}{\displaystyle}
\newcommand{\cal}{\mathcal}

\newcommand{\sub}{\subseteq}
\newcommand{\led}[4]{$\xymatrix@1{{#1\,} \ar[q-1]^-{#2}_-{#3}& {\,#4}}$}
\newcommand{\wed}[4]{$\xymatrix@C=27pt@1{{#1\,} \ar@{->>}[q-1]^-{#2}_-{#3}& {\,#4}}$}

\begin{document}

\title[meet-irreducible discrete copulas]
{Matrix representation of meet-irreducible discrete copulas}
\author{Masato Kobayashi\\
}
\thanks{e-mail: \texttt{kobayashi@math.titech.ac.jp}}
\thanks{Graduate School of Science and Engineering Department of Mathematics\\
Saitama University, 255 Shimo-Okubo, Saitama 338-8570, Japan.
}
\thanks{M. Kobayashi, Matrix representation of meet-irreducible discrete copulas, Fuzzy Sets and Systems, vol. 240 (2014), 117-130.}
\date{\today}
\address{Graduate School of Science and Engineering\\
Department of Mathematics\\
Saitama University,
255 Shimo-Okubo, Saitama 338-8570, Japan.}
\email{kobayashi@math.titech.ac.jp}

\begin{abstract}
Following Aguil\'{o}-Su\~{n}er-Torrens (2008), Koles\'{a}rov\'{a}-Mesiar-Mordelov\'{a}-Sempi (2006)
and Mayor-Su\~{n}er-Torrens (2005), we continue to develop a theory of matrix representation for discrete copulas. To be more precise, we give characterizations of meet-irreducible discrete copulas from an order-theoretical aspect: we show that the set of all irreducible discrete copulas is a lattice in analogy with Nelsen and \'{U}beda-Flores (2005). Moreover, we clarify its lattice structure related to Kendall's $\tau$ and Spearman's $\rho$ borrowing ideas from Coxeter groups.
\end{abstract}

\maketitle

\section{Introduction}

\subsection*{Copulas}

The theory of copulas has been of fundamental importance in probability and statistics.
It dates back to Sklar's Theorem (1959) \cite{sklar}:
\bft{Let $X, Y$ be two random variables with marginal distribution functions $F$ and $G$.
Then there exists a copula $C:[0, 1]^2\to [0, 1]$ such that we can express the joint distribution $H$ of $X$ and $Y$ as $H(x, y)=C(F(x), G(y))$.
}\eft
Since then, we have continued to develop this theory extensively. 
Quasi-copulas, a more general concept, recently appeared in Alsina-Nelsen-Schweizer (1993) \cite{ans}. These days this idea has wide applications in other areas such as Fuzzy logic and Quantitative finance. 

\subsection*{Motivation}
In this article, we focus on a certain class of copulas, \emph{discrete couplas}. Why do we study this class? Here we list some results to see its importance:
\begin{itemize}
\item It has a lattice structure. This is analogous to Nelsen and \'{U}beda-Flores \cite{nu} that all non-discrete copulas have a lattice structure. 
\item \cite[Theorem 2.4]{qs} Every quasi-copula is a certain limit of discrete ones.
\item Discrete quasi-copulas contain rich mathematical structures of not only a matrix but also a group, a lattice and a vector space. In particular, a certain partial order fits into a framework of the classic topic on Kendall's $\tau$ and Spearman's $\rho$ such as Daniels \cite{daniels}, Durbin-Stuart \cite{ds}, Kruskal \cite{kruskal}, Lehmann \cite{lehmann} and Okamoto-Yanagitmoto \cite{oy}.
\item We can make the best use of matrices; as in the title, our study continues recent (2000s) developments on matrix representation of discrete copulas.
\end{itemize}

In mathematics, it is a common idea to understand general objects in terms of ``smaller" ones;
for example, each natural number is a product of prime numbers; each element in a vector space is a linear combination of vectors in a basis. Furthermore, such expressions are often unique.


\bq{What about discrete copulas?  Is this sort of argument possible?
}\eq

The answer is yes. As mentioned above, discrete copulas form a lattice.
We then come to the fundamental fact in the lattice theory: in a finite distributive lattice, each non-maximal element is the meet of meet-irreducible elements. However, we could not find any references on discrete copulas from this aspect in spite of its importance. 
%


\subsection*{State of the art (Matrix representation)}

\begin{figure}
\caption{recent developments on discrete quasi-coplulas}
 \[ \xymatrix@=5mm{
*++<1pt>[F]\txt{This article \\
meet-irreducible copulas\\
$\updownarrow$\\MM
}\ar@{->}_*\txt{Coxeter group}[dd]
\ar@{->}^*\txt{Dedekind-MacNeille\\ completion}[rrdd]
&&\\
&&&\\
*++<1pt>[F]\txt{Mayor-Su\~{n}er-Torrens (2005)\\irreducible copulas\\ $\updownarrow$\\
 PM}\ar@{->}_*\txt{Convex closure}[dd]
& *{\subseteq}
&*++<1pt>[F]\txt{Aguil\'{o}-Su\~{n}er-Torrens (2008)\\irreducible quasi-copulas\\  $\updownarrow$\\ ASM}\ar@{->}^*\txt{Convex closure}[dd]\\
&&&&
\\
*++<1pt>[F]\txt{ Koles\'{a}rov\'{a}-Mesiar-Mordelov\'{a}-Sempi (2006)\\copulas\\  $\updownarrow$\\ BM}&
*{\subseteq} &*++<1pt>[F]\txt{Aguil\'{o}-Su\~{n}er-Torrens (2008)\\
quasi-copulas\\  $\updownarrow$\\ GBM}
}\]
\label{gbb}
\end{figure}
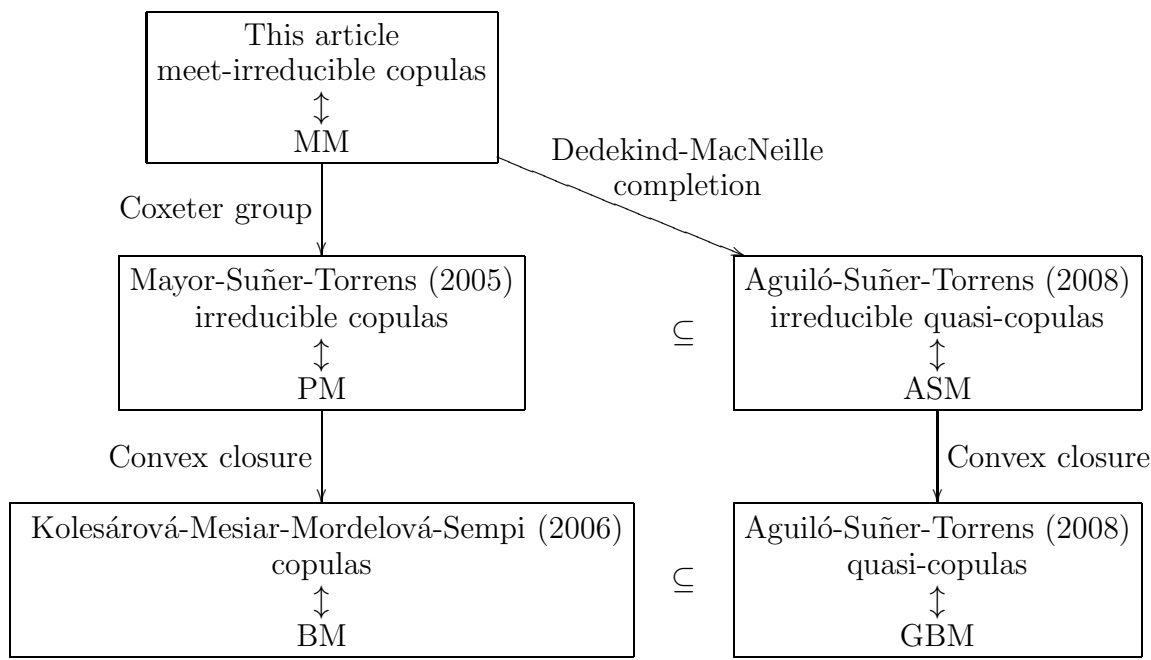

Figure \ref{gbb} shows some state of the art on matrix representation of discrete copulas; there, we see the correspondences between five kinds of discrete quasi-copulas and square matrices: MM means meet-irreducible matrices, PM permutation matrices, ASM alternating sign matrices, BM bistochastic matrices and GBM generalized bistochastic matrices; see also \cite{ast2} for non-square variants. Although we deal with only the first three classes, it is easy to write down several consequences for BM and GBM; details will appear in a subsequent publication.
Our specific goal is to give explicit descriptions of \emph{meet-irreducible copulas} with matrix representations.
For this purpose, we ``borrow" some ideas from Coxeter groups and Bruhat order.



\emph{Coxeter groups} play a significant role in algebra, combinatorics and geometry (initiated by H.S.M. Coxeter and afterward developed by Bourbaki around 1960). 
Symmetric groups are indeed type A Coxeter groups equipped with a certain partial order, called \emph{Bruhat order}. 
%
The key idea with a connection to copulas is the following:

\begin{fact}
Concordance order on irreducible discrete copulas on $\{0, 1, \dots, n\}$ is isomorphic to reverse Bruhat order on the symmetric group $S_n$.
\end{fact}
This order plays a key role for studying a lattice structure of discrete copulas together with matrix representation, as we shall see.
\subsection*{Outline of the paper}
In Section 2, we first set up definitions of discrete copulas and quasi-copulas.
After that, following \cite{ast,kmms,mst}, we recall correspondences between copulas and matrices. Then we go into their poset structures. We show in Theorem \ref{latt} that the set of discrete irreducible quasi-copulas is a lattice in analogy with \cite{nu}.
In Section 3, we give explicit characterizations of meet-irreducible copulas in terms of matrix entries. We then observe some consequences on algebaic properties and Kendall's $\tau$.
Finally, Theorem \ref{mth1} clarifies the relation between entries of discrete copulas and the lattice structure.
In Appendix, we give some basic terminology on posets and Coxeter groups.

\section{Copulas}


\subsection*{Discrete quasi-copulas}
To begin with, let us fix a positive integer $n$.
Below, we treat only a special class of quasi-copulas, \emph{discrete} quasi-copulas on $L_n=\{0, 1, \dots, n\}$ (not on $I_n=\{0, \fr{1}{n}, \dots, 1\}$); we even omit the word ``discrete" whenever no confusion arises.
 For definitions and characterizations of general (quasi-)copulas, see Nelsen's book \cite{nelsen} as well as several recent papers \cite{ans,gqrs,ru}.\\
Let $L$ denote the closed interval $[0, n]$ in $\rr$. 



 \bd{A binary operation $C:L_n\times L_n\to L$ is a (two-dimensional) \emph{discrete copula} if it satisfies all of the following:
\begin{description}
\item[C1] $C(i, 0)=C(0, i)=0$ for all $i\in L_n$.
\item[C2] $C(i, n)=C(n, i)=i$ for all $i\in L_n$.
\item[C3] $C(i, j)+C(i', j')\ge C(i, j')+C(i', j)$ for all $i\le i'$ and $j\le j'$ in $L_n$. (2-increasing)
\end{description}
}\ed
\bd{Say \label{qc}$Q:L_n\times L_n\to L$ is a \emph{discrete quasi-copula} if it satisfies the following:
\begin{description}
\item[Q1] $Q(i, 0)=Q(0, i)=0$ and $Q(i, n)=Q(n, i)=i$ for all $i\in L_n$.
\item[Q2] $Q$ is non-decreasing in each component.
\item[Q3] $Q(i, j)+Q(i', j')\ge Q(i, j')+Q(i', j)$ for all $i\le i'$ and $j\le j'$ in $L_n$ such that at least one of $\{i, i', j, j'\}$ is $0$ or $n$. (2-increasing on the boundary)
\end{description}
}\ed
A discrete quasi-copula which is not a copula is called a \emph{proper} quasi-copula.\\
Thanks to \textbf{C2} and \textbf{Q1}, the range of such operations always contains $L_n$. 
\bd{A discrete copula is \emph{irreducible} if its range is exactly $L_n$.
Similarly, a discrete quasi-copula is \emph{irreducible} if its range is $L_n$.}\ed

Denote by $\cal{P}_n$ ($\mathcal{Q}_n$) the set of all irreducible discrete (quasi-)copulas on $L_n$.
These are our main objects in the sequel.


\brk{In the literature, it is more common to define these operations on $I_n=\{0, \fr{1}{n}, \dots, 1\}$ instead of $L_n$. However, there is no significant difference between them. It is just up to a scale shift; see \cite[Remark 4]{ast}. 
A merit of working on $L_n$ is that we can make the best use of matrices to represent quasi-copulas;
identify a quasi-copula $Q$ with the matrix whose $(i, j)$-entry is $Q(i, j)$. Although some authors prefer to``$xy$-axis notation", we stick to this matrix representation throughout.
Note: for simplicity, we omit $0$ values in the zero-th row and column. 
For example, the matrix $\left[\begin{array}{ccc}  0& 1  &1   \\ 0 & 1  &2   \\ 1 & 2  &3  \end{array}\right]$ represents the quasi-copula $Q$ such that $Q(1, 1)=0, Q(1, 2)=Q(1, 3)=1$ and so on. 


}\erk

\subsection*{Permutation matrices (PM)}
We use notation $[n]=\{1, 2, \dots, n\}$. Let $A=(a_{ij})$ be an $n\times n$ matrix.

\bd{Let $A$ be a \emph{permutation matrix} (PM): for all $(i, j)\in [n]\times [n]$, we have $a_{ij}\in \{0, 1\}$, 
\[\sum_{k=1}^j a_{ik}\in \{0, 1\}, \sum_{k=1}^i a_{kj}\in \{0, 1\}
\mbox{ and }
\sum_{k=1}^n a_{ik}=\sum_{k=1}^n a_{kj}=1.
\] }\ed
In other words, a 1 appears exactly once in every row and column of $A$ and all other entries are 0.\\[.1in]
For example, $
\begin{bmatrix}
0&1&0\\
0&0&1\\
1&0&0
\end{bmatrix}
$ is a PM.
The permutation 231 tells positions of 1's in each row.
Thus, it is often convenient to refer to such matrices as $A(231)$. 
Obviously, PMs have a group structure; the matrix multiplication corresponds to the composition of permutations; inverse matrices therefore correspond to inverse permutations.
Under this identification, we sometimes write $S_n$ (the symmetric group on $[n]$) to mean the set of PMs of size $n$, by slight abuse of language.

\bft{\cite[Proposition 6]{mst}
A binary operation $C$ on $L_n$ is an irreducible discrete copula if and only if 
there exists a unique PM $A=(a_{ij})$ such that 
\[C(r, s)=\begin{cases}
0& \mbox{if $r=0$ or $s=0$},\\
\di\sum_{\substack{i\le r\\ j\le s}}a_{ij}&\mbox{otherwise},
\end{cases}\]
for all $(r, s)\in L_n\times L_n$
and $a_{ij}=C(i, j)+C(i-1, j-1)-C(i, j-1)-C(i-1, j)$
for all $(i, j)\in [n]\times [n]$.
}\eft

This correspondence $C\leftrightarrow A$ is a bijection: $\cal{P}_n\cong S_n$ as sets.
Later, we improve this to an isomorphism of posets.
\subsection*{Alternating sign matrices (ASM)}
\bd{Let $A$ be an \emph{alternating sign matrix} (ASM): for 
all $(i, j)\in [n]\times [n]$, we have $a_{ij}\in \{-1, 0, 1\}$, 
\[\sum_{k=1}^j a_{ik}\in \{0, 1\}, \sum_{k=1}^i a_{kj}\in \{0, 1\}
\mbox{ and }
\sum_{k=1}^n a_{ik}=\sum_{k=1}^n a_{kj}=1.
\]
Denote by $\tn{ASM}_n$ the set of all ASMs of size $n$.
}\ed
Observe that every PM is an ASM. Say an ASM is \emph{proper} if it is not a PM.\\
For example, 
$
\begin{bmatrix}
0&1&0\\
1&-1&1\\
0&1&0
\end{bmatrix}
$ is a proper ASM. Note that this ASM does not have the matrix inverse.
As we shall see, instead of losing a group structure, ASMs behave order-theoretically much better.


\bft{\cite[Proposition 8]{ast}
A binary operation $Q$ on $L_n$ is an irreducible discrete quasi-copula if and only if 
there exists a unique ASM $A=(a_{ij})$ such that 
\[Q(r, s)=\begin{cases}
0& \mbox{if $r=0$ or $s=0$},\\
\di\sum_{\substack{i\le r\\j\le s}}a_{ij}&\mbox{otherwise},
\end{cases}\]
for all $(r, s)\in L_n\times L_n$ and 
$a_{ij}=Q(i, j)+Q(i-1, j-1)-Q(i, j-1)-Q(i-1, j)$
for all $(i, j)\in [n]\times [n]$.
}\eft

Again, this correspondence $Q\leftrightarrow A$ is a bijection: $\cal{Q}_n \cong \tn{ASM}_n $ as sets. Later, we improve this to an isomorphism of posets.

\brk{On the one hand, the number of PMs of size $n$ is $n!$.
On the other hand, the numbers of ASM's are \emph{Robbins numbers}:
\[A_n=\prod_{i=0}^{n-1} \frac{(3i+1)\,!}{(n+i)\,!}.\]
This sequence goes as $1, 2, 7, 42, 429, 7436, 218348, 10850216, \dots$; see Robbins' article \cite{robbins} on the long history for discovery of this formula.
Evidently, it is difficult to deal with such a large number of matrices; we want to find a smaller number of ``nice" matrices to study general quasi-copulas.
This desire leads us to another class of matrices (and copulas) which we call \emph{meet-irreducible}. It will be clear that this class is appropriate when we investigate a certain partial order on $\cal{Q}_n$.
Thus we have to mention that order next.
}\erk





\subsection*{Concordance and reverse Bruhat orders}

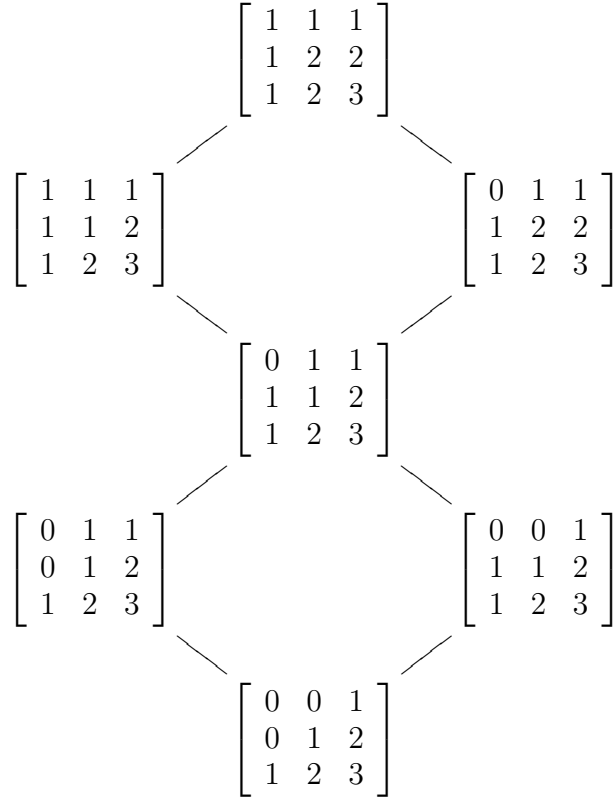
\begin{figure}
\caption{$(\mathcal{Q}_3, \le)$}
\label{q3}
\[
\xymatrix@=5mm{
&*+{\left[\begin{array}{ccc}1  &1   &1   \\ 1 & 2  & 2  \\1  & 2  &3  \end{array}\right]}
\ar@{-}[dl]\ar@{-}[dr]
&\\
*+{\left[\begin{array}{ccc}1  &1   &1   \\1  &1   &2   \\ 1 & 2  &3  \end{array}\right]}
\ar@{-}[dr]
&&*+{\left[\begin{array}{ccc}0 &1   &1   \\ 1 &2  &2   \\1  & 2  &3  \end{array}\right]}\ar@{-}[dl]
\\
&*+{\left[\begin{array}{ccc} 0 &1   &1   \\1  &1   &2   \\1 &2   &3  \end{array}\right]}
\ar@{-}[dl]\ar@{-}[dr]
\\
*+{\left[\begin{array}{ccc} 0 & 1  &1   \\ 0 &1  &2   \\1  &2   &3  \end{array}\right]}
\ar@{-}[dr]
&&*+{\left[\begin{array}{ccc}  0& 0  &1   \\ 1 & 1  &2   \\1  &2   &3  \end{array}\right]}
\ar@{-}[dl]
\\
&*+{\left[\begin{array}{ccc}  0& 0  &1   \\0  &1   &2   \\ 1 &2   &3  \end{array}\right]}
&
}\]
\end{figure}

\bd{Let $P, Q\in \cal{Q}_n$.
Define the (discrete) \emph{concordance order} $P\le Q$ if $P(i, j)\le Q(i, j)$ for all $i, j \in L_n$.
}\ed
Note that for $i=0$ or $j=0$, the equality $P(i, j)=Q(i, j)$ always holds because both are $0$. Hence we do not care about these entries.\\
Figure \ref{q3} illustrates this ``entrywise" order for $n=3$;
from the bottom to top, entries get larger and larger one by one, as we see.

\brk{Let $|Q|$ denote a sum of all entries of $Q\in \cal{Q}_n$.
By definition of $<$ above, 
if $P\lhd Q$ (a covering relation), then $|Q|-|P|$ is a positive integer; here we are not saying that $|Q|-|P|$ is $1$, though. This is indeed the case as shown in Theorem \ref{mth1}. In fact, $Q$ is a \emph{graded} poset (Appendix \ref{apd}).
}\erk

We now come to the poset isomorphism as mentioned before. On the one hand, the subset $\cal{P}_n$ ($\sub \cal{Q}_n$) naturally has the induced suborder.
On the other hand, $S_n$ has the reverse Bruhat order $\le'$ (Appendix \ref{apd}).

\begin{thm}

The posets $(S_n, \le')$ and $(\cal{P}_n, \le)$ are isomorphic.
\end{thm}


\bpf{Let $w\in S_n$ and consider the corresponding PM, say $A(w)=(a_{ij})$. On the PM side, $w(r, s)=|\{i\mid i\le r \mbox{ and }w(i)\le s\}|$ counts the number of
 $1$'s in the northwest part of an $(r, s)$-entry (the border inclusive) of $A(w)$.
On the copula side, $Q(r, s)$ is the sum $\sum_{i\le r, j\le s} a_{ij}$.
But this quantity coincides with $w(r, s)$ because each $a_{ij}$ is either $0$ or $1$. Therefore, entrywise orders determined by $\{w(r, s)\}$ and $\{Q(r, s)\}$ are same.
}\epf

In what follows, we use the two orders interchangeably.


Next we go into details of a poset structure of $\cal{Q}_n$;
below Theorem \ref{latt} asserts that $\cal{Q}_n$ is a lattice. Here we need two lemmas:

\bl{\label{q1}
Let $Q:L_n\times L_n\to L$. Then the following are equivalent:
\be{
\item $Q\in \cal{Q}_n$.
\item $Q$ satisfies \textbf{Q1} and \label{x2}
moreover
\[Q(i, j)-Q(i-1, j)\in \{0, 1\} \mbox{ and } Q(i, j)-Q(i, j-1)\in \{0, 1\}\]
for all $i, j\in [n]$.
}\ee
}\el

\bpf{This is a consequence of \cite[Proposition 2 and Remark 4]{ast} .
}\epf

\bl{Let $P, Q\in \cal{Q}_n$.
Define a matrix $R$ by $R(i, j)=\min\{ P(i, j), Q(i, j)\}$.
Then $R\in\cal{Q}_n$ and $R=P\we Q$, that is, $R$ is the meet of $P$ and $Q$ in the poset $(\mathcal{Q}_n, \le)$.
}\el
\bpf{
We need to show that (a): $R$ is indeed a quasi-copula and 
(b): $R$ is the maximum element of the set $\{S\in \cal{Q}_n\mid S\le P \mbox{ and } S\le Q\}$.\\
\indent (a): For convenience, we use Lemma \ref{q1} instead of verifying Definition 
\ref{qc}.
It is easy to check \textbf{Q1}: $R(i, 0)=\min \{P(i, 0), Q(i, 0)\}=0$ ($=R(0, i)$) and $R(i, n)=\min \{P(i, n), Q(i, n)\}=i$ ($=R(n, i)$). Next, fix $i, j \in [n]$. 
To show that $R(i, j)-R(i-1, j)\in \{0, 1\}$, let $\al=P(i-1, j), \beta=P(i, j)$, $\gam=Q(i-1, j)$ and $\delta=Q(i, j)$ for simplicity.
Then 
$\beta-\al\in\{0, 1\}$ and $\delta-\gam\in \{0, 1\}$ because $P$ and $Q$ satisfy Lemma \ref{q1} (\ref{x2}).
(i) Suppose $\delta=\beta$.
Then $R(i, j)-R(i-1, j)=\beta-\min\{\al, \gam\}\in \{0, 1\}$.
(ii) Suppose $\delta-\beta\ge 1$. Then 
$\alpha\le \beta\le \delta-1\le \gam$ so that 
$R(i, j)-R(i-1, j)=\min\{\beta, \delta\}-\min\{\al, \gam\}=\beta-\al\in \{0, 1\}$.
(iii) For the case $\beta-\delta\ge 1$, we can mimic the proof above. 
Finally, we can show $R(i, j)-R(i, j-1)\in \{0, 1\}$ in the same way.\\
\indent (b): Suppose $S\le P \mbox{ and } S\le Q$.
It follows that $S(i, j)\le P(i, j)$ and $S(i, j)\le Q(i, j)$ for all $i, j$.
Hence $S(i, j)\le \min \{P(i, j), Q(i, j)\}$ for all $i, j$. In other words, $S\le R$.
}\epf

\bt{For $n\ge 1$, $\cal{Q}_n$ is a lattice. \label{latt}
For $n\ge 3$, $\cal{P}_n$ is not a lattice.
For $n\ge 4$, $\cal{Q}_n\setminus \cal{P}_n$ is not a lattice.
}\et


\bpf{We just proved that $P\we Q$ exists and belongs to $\cal{Q}_n$. With the order-dual argument, $P\vee Q$ also exists and belongs to $\cal{Q}_n$.
As a result, $\cal{Q}_n$ is a lattice.
To show that $\cal{P}_n$ is not a lattice for $n\ge 3$, we first consider $P, Q\in \cal{P}_3$ corresponding to $v=132$ and $w=213$. The matrix $P \we Q$ is a proper quasi-copula (the middle matrix in Figure \ref{q3}). For $n\ge 4$, consider the embedding of these permutations: $v'=13245\cdots n$ and $w'=21345\cdots n$. Then discuss $v' \we w'$ in $S_n$.
It remains to show that $\cal{Q}_n\setminus \cal{P}_n$ is not a lattice for $n\ge 4$.
Here we claim that there exist $P, Q \in \cal{Q}_4\setminus \cal{P}_4 $ such that $P\we Q \not\in \cal{Q}_4\setminus\cal{P}_4$. 
Let  $A=\left[\begin{array}{cccc}0  & 0  & 0  & 1  \\ 0 & 1  & 0  & 0  \\  1&  -1 & 1  & 0  \\  0& 1  &0   & 0 \end{array}\right]$
and $B=\left[\begin{array}{cccc}0  & 0  & 1  & 0  \\0  & 1  & -1  & 1  \\0  & 0  & 1  & 0  \\ 1 & 0  & 0  & 0  \end{array}\right]$ be proper ASMs. They correspond to proper quasi-copulas  
$P=
\left[\begin{array}{cccc}0  & 0  & 0  & 1  \\ 0 & 1  & 1  & 2  \\  1&  1 & 2  & 3  \\  1& 2  &3   & 4 \end{array}\right]$
 and 
$Q=\left[\begin{array}{cccc}0  & 0  & 1  & 1  \\0  & 1  & 1  & 2  \\0  & 1  & 2  & 3  \\ 1 & 2  & 3  & 4  \end{array}\right]$.
Observe that 
$P\we Q=\left[\begin{array}{cccc} 0 & 0  &  0 &  1 \\ 0 & 1  & 1  & 2  \\ 0 & 1  & 2  & 3  \\ 1 & 2  & 3  &4  \end{array}\right]$
which corresponds to the PM
$\left[\begin{array}{cccc} 0 & 0  &  0 &  1 \\ 0 & 1  & 0  & 0  \\ 0 & 0  & 1  & 0  \\ 1 & 0  & 0  &0  \end{array}\right]$.
For $n\ge 5$, again, consider the natural embedding.
}\epf

\begin{cor}\label{fin}
As the Dedekind-MacNeille completion of $(\cal{P}_n, \le)\cong (S_n, \le')$, we obtain $(\cal{Q}_n, \le) \cong (\tn{ASM}_n, \le')$.
\end{cor}

In fact, $\cal{Q}_n$ is a finite distributive lattice (Appendix \ref{apd}).

\brk{Nelsen and \'{U}beda-Flores \cite{nu} showed that the set of all non-discrete quasi-copulas $\cal{Q}$ is the Dedekind-MacNeille completion of copulas $\cal{C}$ and moreover $\cal{Q}$ is a lattice whereas neither $\cal{C}$ nor $\cal{Q}\setminus \cal{C}$ is a lattice.
Hence we may understand that our theorem and corollary above is a discretization of their results.
}\erk

Consequently, $\cal{Q}_n$ has the minimum and maximum elements (denoted by $W_n$ and $M_n$), \emph{Fr\'{e}chet-Hoeffding bounds}: $W_n(i, j)=\max\{i+j-n, 0\}$ and $M_n(i, j)=\min\{i, j\}$.
In fact, $M_n$ and $W_n$ are  elements of $\cal{P}_n$ corresponding to the reverse and identity permutations, respectively. 
Despite such simple definitions, these copulas play a key role in the proof of Theorem \ref{mth1}.


\section{Meet-irreducible discrete copulas}

\subsection*{Meet-irreducibility}\label{mirr}
We now introduce our main idea, meet-irreducibility. Let us introduce its definition in a rather general setting.
\bd{
Let $x, y, z$ be elements of a finite poset $(P, \le)$. Define $z$ to be \emph{meet-irreducible} in $P$ if $z$ is not the maximum element of $P$ and whenever $z=x \wedge y$ then $z=x$ or $z=y$. Denote by $M(P)$ the set of all meet-irreducible elements in $P$.
}\ed

It follows that each non-maximal $z\in P$ can be written as $z=x_1\wedge \dots \wedge x_k$
for some anti-chain $\{x_1, \dots, x_k\}$ in $M(P)$ (take all minimal elements of $\{x\in M(P) \mid z\le x\}$); see Reading \cite[Section 2]{reading4} for details of the lattice theory. Thus, to know a poset structure of $P$, it is essential to understand such elements.

\bq{
What are $M(\mathcal{P}_n)$ and $M(\mathcal{Q}_n)$?
}\eq

Indeed $M(\cal{P}_n)=M(\cal{Q}_n)$; see Appendix \ref{apd}.
In particular, every meet-irreducible discrete irreducible quasi-copula is necessarily a copula.
\brk{We see an obvious conflict of terms ``meet-irreducible" and ``irreducible".
We will just say ``meet-irreducible discrete copulas" whenever no confusion arises.
}\erk

Now we want to find a necessary and sufficient condition for $Q=Q(A)\in \cal{Q}_n$ to be meet-irreducible in terms of its entries. 
To describe the matrix correspondence $Q\leftrightarrow A$ more explicitly, let us introduce this definition.

\bd{A pair $(i, j)$ is a \emph{positive position} of $Q$ if $a_{ij}>0$ (i.e., $a_{ij}=1$); $(i, j)$ is a \emph{negative position} of $Q$ if $a_{ij}<0$ (i.e., $a_{ij}=-1$).
}\ed

A positive position indicates validity of a strict inequality of the 2-increasing condition for the square with $(i-1, j-1)$, $(i, j-1)$, $(i-1, j)$ and $(i, j)$-entries. Similarly, a negative position indicates a failure of the 2-increasing condition.

\bex{
Let $Q=\left[\begin{array}{ccccc}
0  & 0  & 0  & 1  & 1  \\
 0 & 0  & 1  & 1  & 2  \\
 0   & 1  & 1  & 2  & 3  \\ 
 1   & 1  & 2  & 3  & 4  \\
 1      & 2  & 3  & 4  & 5 \end{array}\right]``="
\left[\begin{array}{ccccc}
  &   &   & \cci &   \\
  &   &  \cci  & \times  &  \cci \\
    &  \cci  & \times  &  \cci  &   \\ 
  \cci   & \times  &  \cci  &   &   \\
       &  \cci  &   &   &  \end{array}\right]$.
Here \cci\phantom{ }shows some positive and $\times$ some negative positions.

}\eex


\bp{Let $Q\in \cal{Q}_n$.\label{mtp}
Then $Q\in M(\cal{Q}_n)$ if and only if 
its positive positions are as below and it does not contain any negative positions.
\[
\left[\begin{array}{ccccccccccccc}
\cci  &   &   &   &   &   &   &   &   \\
  &  \ddots &   &   &   &   &   &   &   \\ 
   &   &  \cci &   &   &   &   &   &   \\  
  &&& &   &   &&  \cci &   &   &   &   &   \\ 
  &&& &   &   &&   &  \ddots &   &   &   &   \\ 
  &&& &   &   &  & &   & \cci  &   &   &   \\ 
   &   &   &   & \cci  &   &   &   &   \\ 
   &   &   &   &   & \ddots  &   &   &   \\  
   &   &   &   &   & & \cci    &   &   \\  
&&&  &  &   &   &   &   &   &  \cci &   &   \\ 
  &&&&   &   &   &   &   &   &   &  \ddots &   \\  
   &&&&  &   &   &   &   &   &   &   &\cci  
     \end{array}\right]
\]
As a result, there do not exist meet-irreducible discrete proper quasi-copulas.
A more precise statement is:
there exist integers $i_1, i_2, i_3, i_4$ and $w\in S_n$ such that 
$i_1+i_2+i_3+i_4=n$,
$i_1\ge 0, i_2, i_3\ge 1, i_4\ge 0$, 
\[w(i)=\begin{cases}
i&\mbox{if }1\le i\le i_1\\
i+i_3&\mbox{if }i_1+1\le i\le i_1+i_2\\
i-i_2&\mbox{if }i_1+i_2+1\le i\le i_1+i_2+i_3\\
i&\mbox{if }i_1+i_2+i_3+1\le i\le n\\
\end{cases}\]
and $A(w)$ is the associated PM for $Q$ (write this permutation as $w=w(i_1, i_2, i_3, i_4)$).
}\ep

\brk{We can rephrase this more intuitively: let $1, 2, \dots, n$ be an increasing sequence. 
Prepare four empty boxes as \fbox{\phantom{1}}\fbox{\phantom{1}}\fbox{\phantom{1}}\fbox{\phantom{1}}\,.
Then put each number of $\{1, 2, \cdots, n\}$ into these boxes keeping the total order, left to right. Our rule is that we allow the first and fourth boxes to be empty while the second and third ones cannot be empty. Then interchange the second and third boxes together with numbers.  Resulting sequences are precisely one-line expressions of meet-irreducible permutations. Examples of $M(S_8)$ are \[\fbox{1} \fbox{34567} \fbox{2} \fbox{8}\,, 
\fbox{1} \fbox{5678} \fbox{234} \fbox{\phantom{9}}\,, 
\fbox{\phantom{0}} \fbox{45} \fbox{123} \fbox{678} \mbox{ and } 
\fbox{\phantom{0}} \fbox{678} \fbox{12345} \fbox{\phantom{1}}\,.\]
}\erk

In this way, it is easy to construct meet-irreducible permutations and hence meet-irreducible copulas.
We postpone the proof to the next subsection.

\begin{ob}\label{there}
For each non-maximal $Q\in \mathcal{Q}_n$, there exist a unique anti-chain of meet-irreducible copulas $R_1, \dots, R_k$ such that $Q=R_1\we \cdots \we R_k$.
\end{ob}

For example, let $Q=\left[\begin{array}{cccc}0 & 0 & 1 & 1 \\0 & 1 & 1 & 2 \\0 & 1 & 2 & 3 \\1 & 2 & 3 & 4\end{array}\right]$. Then

\begin{align*}
Q&=\left[\begin{array}{cccc}0 & 1 & 1 & 1 \\0 & 1 & 2 & 2 \\0 & 1 & 2 & 3 \\1 & 2 & 3 & 4\end{array}\right]\we
\left[\begin{array}{cccc}0 & 0 & 1 & 1 \\1 & 1 & 2 & 2 \\1 & 2 & 3 & 3 \\1 & 2 & 3 & 4\end{array}\right]\we
\left[\begin{array}{cccc}1 & 1 & 1 & 1 \\1 & 1 & 1 & 2 \\1 & 2 & 2 & 3 \\1 & 2 & 3 & 4\end{array}\right]\\
&=
\left[\begin{array}{cccc} & \cci&  &  \\ &  & \cci &  \\ &  &  & \cci \\\cci &  &  & \end{array}\right]\we
\left[\begin{array}{cccc} &  & \cci &  \\\cci &  &  &  \\ & \cci &  &  \\ &  &  & \cci\end{array}\right]\we
\left[\begin{array}{cccc}\cci &  &  &  \\ &  &  &\cci \\ & \cci &  &  \\ &  & \cci & \end{array}\right].
\end{align*}



\subsection*{Meet-irreducible matrices (MM)}

As an analogy of PM and ASM, it is natural to introduce the following (see back Figure \ref{gbb}):
\bd{Let $A$ be a \emph{meet-irreducible matrix} (MM): it is a PM and the associated pemutation is meet-irreducible. 
}\ed

We now bring ideas from Coxeter group theory; there is a simple characterization of meet-irreducible permutations.
For $w\in S_n$, let
$D_L(w)=\{  i \in [n-1] \mid w^{-1}(i)>w^{-1}(i+1) \}$
and  $D_R(w)=D_L(w^{-1})$.
Call these sets \emph{left} and \emph{right descents} of $w$.


\bft{\cite[Sections 7 and 8]{reading4} The following are equivalent:\label{bg}
\be{\item $w$ is meet-irreducible in $(S_n, \le')$. 
\item $|D_L(w)|=|D_R(w)|=1$.\label{bg2}
}\ee
}\eft
In particular, the identity permutation is \emph{not} meet-irreducible.

\brk{
Apart from Bruhat and reverse Bruaht orders, it is common to call permutations satisfying Fact \ref{bg} (\ref{bg2}) \emph{bigrassmannian} in the Coxeter group context.
}\erk



\bp{Let $w\in S_n$. Then the following are equivalent:
\be{\item there exist integers $(i_1, i_2, i_3, i_4)$ as stated in Proposition \ref{mtp}.\label{w1}
\item $|D_L(w)|=|D_R(w)|=1$.\label{w2}
}\ee
}\ep

\bpf{(\ref{w1}) $\Longrightarrow$ (\ref{w2}): Check that $D_L(w)=\{i_1+i_3\}$ and $D_R(w)=\{i_1+i_2\}$.\\
(\ref{w2}) $\Longrightarrow$ (\ref{w1}):
Let $k$ be the unique element of $D_R(w)$.
Thus the image of $w$ splits into two increasing sequences:
$\{w(1)<\cdots<w(k)\}$ and $\{w(k+1)<\cdots<w(n)\}$.
Keeping this in mind, let
\begin{align*}
I&=\{i \mid 1\le i\le k-1 \mbox{ and }w(i+1)-w(i)>1\},\\
J&=\{i \mid k+2\le i\le n \mbox{ and }w(i)-w(i-1)>1\}.
\end{align*}
Some of these sets may be empty. We claim that $|I|, |J|\le 1$.\\
\emph{Proof of Claim.}
Suppose, toward a contradiction, $|I|\ge 2$, say $a, b\in I$ and $a<b$.
Then $w(a+1)-w(a)>1$ and $w(b+1)-w(b)>1$.
Since $w(1)<\dots <w(k)$, we must have $w^{-1}(w(a+1)-1)\not\in [k]$ and $w^{-1}(w(b+1)-1)\not\in [k]$. Therefore $w(a+1)-1, w(b+1)-1 \in D_L(w)$ with $w(a+1)-1\ne w(b+1)-1$, a contradiction.
Similarly, suppose $|J|\ge 2$, say $c, d\in J$ and $c<d$.
Then $w(c)-w(c-1)>1$ and $w(d)-w(d-1)>1$.
Since $w(k+1)<\cdots <w(n)$, we must have $w^{-1}(w(c-1)+1)\in [k]$ and $w^{-1}(w(d-1)+1)\in [k]$. Therefore $w(c-1), w(d-1)$ $\in D_L(w)$, with $w(c-1)\ne w(d-1)$, a contradiction$. \, \blacksquare$\\
Now we determine $i_1, i_2, i_3$ and $i_4$. If $I=\varnothing$, then $i_1:=0$; otherwise, say $I=\{a\}$, let $i_1:=a$.
Similaly, if $J=\varnothing$, then $i_4:=0$; otherwise, say $J=\{c\}$, let $i_4:=n-c$.
Finally, set $i_2:=k-i_1$ and $i_3:=n-i_1-i_2-i_4$ (in any case).
By construction, each of sequences $\{w(1)< \dots <w(i_1)\}, \{w(i_1+1)< \dots <w(i_1+i_2)\}$, 
$\{w(i_1+i_2+1)< \dots <w(i_1+i_2+i_3)\}$ and $\{w(i_1+i_2+i_3+1)< \dots <w(n)\}$ 
is increasing one by one, as required.
}\epf


\bpf[Proof of Proposition \ref{mtp}]
{Suppose $Q\in M(\cal{Q}_n)$, say $Q=Q(w)$ and $w\in M(S_n)$.
Then positive positions of $Q$ are positions of $1$'s in $A(w)$ as explained above. 
Since $Q\in M(\cal{Q}_n)= M(\cal{P}_n) \sub \cal{P}_n$, there is no $-1$ entry in $A(Q)$ so that $Q$ does not contain any negative positions. The converse is clear.
}\epf


\subsection*{Commutativity and associativity}\label{c1}
Here we record some consequences of the last subsection (although we do not need them for Theorem \ref{mth1}) for subsequent research.\\
Define $Q\in \cal{Q}_n$ to be \emph{commutative} if $Q(i, j)=Q(j, i)$ for all $i, j$. Define $Q\in \cal{Q}_n$ to be \emph{associative} if $Q(i, Q(j, k))=Q(Q(i, j), k)$ for all $i, j, k$.
It is easy to characterize these algebraic properties for meet-irreducible copulas:

\begin{ob}
Let $Q\in M(\cal{Q}_n)$, say $Q=Q(w)$ with $w=w(i_1, i_2, i_3, i_4)\in M(S_n)$. 
Then $Q(w)$ is commutative $\iff$ $A(w)$ is symmetric $\iff$ $i_2=i_3$.
\end{ob}
Next, before giving a characterization of associativity for meet-irreducible copulas, we need some definitions.

\bd{The \emph{Lukasiewicz matrix} of size $n$ is the PM $A=(a_{ij})$ with $a_{ij}=1$ whenever $i+j=n+1$.
}\ed

Of course, this matrix corresponds to $W_n$ (Fr\'{e}chet-Hoeffding lower bound) as well as the reverse permutation.

\bft{\cite[Proposition 9]{mst} A discrete copula is associative if and only if its associated PM is an ordinal sum of Lukasiewicz matrices.
}\eft

Define $Q\in \cal{P}_n$ to be \emph{Coxeter} if the associated permutation is a Coxeter generator; these are precisely coatoms of $\cal{P}_n$ and $\cal{Q}_n$.
\bp{
Let $Q\in M(\cal{Q}_n)$. Then $Q$ is associative $\iff$ $Q$ is Coxeter.\label{er}
}\ep

\bpf{Let $A(w)\in M(S_n)$ be the associated PM for $Q$.
Suppose $Q$ is associative.
Then $A=A_1\oplus \dots \oplus A_k$, with each $A_i$ Lukasiewicz matrix.
If the size of some $A_i$ is greater than $2$, then $w$ has more than one right descent, a contradiction since $w$ is meet-irreducible.
If all of size of $A_i$ are 1, then $w$ is the identity permutation, a contradiction again.
Hence $\{A_i\}$ contains a size $2$ matrix, say $A_j$, $1\le j\le n-1$. Moreover, such $j$ must be unique for the same reason. It follows that $A$ must be of the form
\[\left[\begin{array}{cccccccc}
1  &   &   &   &   &   &   &   \\ 
 &  \ddots &   &   &   &   &   &   \\
   &   &1   &   &   &   &   &   \\ 
    &   &   &   &1   &   &   &   \\  
    &   &   &   1&   &   &   &   \\
      &   &   &   &   &1   &   &   \\
        &   &   &   &   &   &\ddots   &   \\  
        &   &   &   &   &   &   & 1 \end{array}\right].\]
Thus $w$ is the transposition interchanging $j$ and $j+1$, a Coxeter generator.
The converse is clear.


}\epf




\subsection*{Kendall's $\tau$}\label{c2}
In the course of studying (non-discrete) copulas, \emph{Kendall's} $\tau$ plays a fundamental role. It takes real values in $[-1, 1]$ satisfying many inequalities with other statistics; see Chapter 5 of Nelsen's book \cite{nelsen} for details. Here we consider some similar statistic for meet-irreducible copulas in our discrete setting; we will show that statistic is ``almost positive" for meet-irreducibles.\\
Say $(i, j)$ is an \emph{inversion} of $w$ if $i<j$ and $w(i)>w(j)$.
Let $\ell(w)$ be the number of inversions of $w$. 

\bd{\[\tau(Q(w))=1-\fr{2\ell(w)}{\fr{n(n-1)}{2}}.\]
}\ed

For a real number $\alpha$, let $\lf \alpha \rf$ denote the least integer which does not exceed $\alpha$.

\bl{\label{lm}
\[\max \{\ell(w)\mid w\in M(S_n)\}=\lf n^2/4\rf.\]
}\el

\bpf{Observe that $\ell(w(i_1, i_2, i_3, i_4))=i_2i_3$.
Hence it is enough to find the maximum for $i_2i_3$ under the condition $i_2\ge 1$, $i_3\ge 1$ and $i_1+i_2+i_3+i_4=n$. Equivalently, find maximal area of rectangles with integer width $i_2$ and height $i_3$ summing up to $n$.
Thus $w(0, \lf(n+1)/2\rf, \lf n/2 \rf, 0)$ gives the maximum $i_2i_3=\lf(n+1)/2\rf\lf n/2\rf=\lf n^2/4\rf$.
}\epf


\bp{
If $Q\in M(\cal{Q}_n)$ $(n\ge 2)$, then 
\[-\fr{1}{n-1}\le \tau(Q)\le 1.\]
}\ep

\bpf{The second inequality is clear. We verify the first one.
Say $Q=Q(w), w\in M(S_n)$. Thanks to Lemma \ref{lm}, we have
\[\tau(Q)=1-\fr{2\ell(w)}{\fr{n(n-1)}{2}}\ge 1-\fr{2 (n^2/4)}{\fr{n(n-1)}{2}}=-\fr{1}{n-1}.\]
}\epf


\subsection*{Rank function}\label{mto}
Another important statistic for non-discrete copulas is Spearman's $\rho$.
As discussed in Nelsen's book \cite[Chapter 5]{nelsen}, $\tau$ and $\rho$ are certain integral of copulas.
\bq{As an analogy, what if we take a sum of values of $Q\in \cal{Q}_n$ in our discrete setting? 
}\eq
Let $|Q|=\sum_{i, j}Q(i, j)$ denote the sum of all entries of the matrix $Q\in \mathcal{Q}_n$; we can still ignore zero-th row and column.
Then $|\phantom{Q}|: \cal{Q}_n\to \rr$ is a strictly increasing function.\\

\bp{$|M_n|=\fr{n(n+1)(2n+1)}{6}$.
}\ep
\bpf{Recall that $M_n(i, j)=\min\{i, j\}$.
Write down its entries as
\[\left[\begin{array}{ccccc}
1  &1   &   1&   \cdots&1   \\
1  & 2  & 2  &\cdots   &2   \\ 
1   & 2 &   &   & \vdots  \\ 
 \vdots   & \vdots  &   & n-1  &  n-1 \\ 
  1   & 2  & \cdots  &n-1   &n  \end{array}\right]\]
so that $  |M_n|=1\, (2n-1)+2\,(2n-3)+\dots+ (n-1)\,3+n\,1=\fr{n(n+1)(2n+1)}{6}.
$
}\epf

\bd{$m(Q)=|M_n|-|Q|$.
}\ed


\bl{$m(W_n)=\fr{(n-1)n(n+1)}{6}=|M(\cal{Q}_n)|$ and $m(M_n)=0$.\label{lem}
}\el

\bpf{
Recall that $W_n(i, j)=\max\{i+j-n, 0\}$.
Write down its entries as
\[\left[\begin{array}{cccccc}
  &   &   &&   &1   \\
  &   &   && 1  & 2  \\  
  &   &  1 &&  \vdots & \vdots  \\ 
   &  1 & \cdots  &&  \vdots &  n-1 \\
    1 &  2 &  \cdots && n-1  & n \end{array}\right].\]
Above the anti-diagonal, all blank positions are $0$. Thus 
\[|W_n|=1\,n+2\,(n-1)+\dots +(n-1)\,2+n \,1=\fr{n(n+1)(n+2)}{6}.\]
Consequently, $m(W_n)=\fr{n(n+1)(2n+1)}{6}-\fr{n(n+1)(n+2)}{6}=\fr{(n-1)n(n+1)}{6}$.
Finally, $m(M_n)=0$ is immediate.
}\epf

The following theorem now reveals the relation between entries of discrete quasi-copulas and the (graded) lattice structure.
\begin{thm}\label{mth1}
$m(Q)=|\{R\in M(\cal{Q}_n)\mid Q\le R\}|$.
\end{thm}

\bpf{
The function $m:\cal{Q}_n\to \rr$ is strictly decreasing by definition of the concordance order.
In particular, if $P\lhd Q$ (a covering relation) then $m(P)-m(Q)$ is not only positive but also an integer (hence at least one). With Lemma \ref{lem}, we know that $m(W_n)=|M(\cal{Q}_n)|$ and $m(M_n)=0$. 
Since $\cal{Q}_n$ is a graded poset of rank $|M(\cal{Q}_n)|$, every maximal chain (from $W_n$ up to $M_n$) has the length exactly $|M(\cal{Q}_n)|$.
Hence for each covering relation $P\lhd Q$, a positive integer $m(P)-m(Q)$ must be $1$. This shows that $m$ coincides with the function $m_{\cal{Q}_n}$ for the distributive lattice $\cal{Q}_n$ (Appendix \ref{apd}).
}\epf
\bex{Let $Q=Q(231)=\left[\begin{array}{ccc}0  &1   &1   \\ 0 & 1  & 2  \\1  & 2  &3  \end{array}\right]$. Then 
$m(Q)=
\left|\begin{array}{ccc}1  &1   &1   \\ 1 & 2  & 2  \\1  & 2  &3  \end{array}\right|
-\left|\begin{array}{ccc}0  &1   &1   \\ 0 & 1  & 2  \\1  & 2  &3  \end{array}\right|
=3$.}\eex


\brk{
If $Q$ is a copula, say $Q=Q(w)$ and $w\in S_n$, then $m(Q(w))$ is equal to 
\[\beta(w):=\sum_{\substack{i<j\\w(i)>w(j)}}(w(i)-w(j)).\] 
See \cite[Theorem]{koba}.
In the example above, $\beta(231)=(2-1)+(3-1)=3$. 
This is a useful formula to compute $m(Q)$.
}\erk

\section{Conclusion}
We have studied discrete copulas from a lattice-theoretic point of view as a continuation of recent work on matrix representations \cite{ast,kmms,mst,nu}. The main idea was to introduce a new class of copulas, \emph{meet-irreducible} copulas. Then we showed its characterization in terms of matrix entries. 
This method clarified lattice structures of ASMs as well as discrete quasi-couplas. We also observed some consequences from algebraic and enumerative aspects such as commutativity and Kendall's $\tau$. In this way, discrete copulas have rich mathematical structures.\\
\indent We end with some ideas for our future research.
\begin{itemize}
\item Observtation \ref{there} guarantees (theoretically) the existence of a decomposition of a given quasi-copula into the meet of meet-irreducible ones. How can we find such a decomposition?
\item It should be possible to develop similar ideas for copulas for \emph{non-square} matrices as studied in \cite{ast2}. For example, it makes sense to speak of the concordance order for such matrices. Study this order in details.
\end{itemize}

\noindent\textbf{Acknowledgments}\\
The author thanks the anonymous referees as well as the editor. Their useful comments and suggestions improved the manuscript.


\appendix

\section{Poset and Coxeter group}
\label{apd}
In this appendix, we recall some definitions and provide useful facts.
Reading \cite{reading4} contains most of these. 

\subsection*{Poset} \label{ps}

Let $(P, \le)$ be a finite poset and $x, y\in P$. 
Say $y$ \emph{covers} $x$ (write $x\lhd y$) if $x<y$ and $\{z\in P \mid x<z<y\}=\varnothing$.
Say $X=\{x_1, \dots, x_n\} \sub P$ is a \emph{chain} if 
whenever $x_i\ne x_j$, then $x_i<x_j$ or $x_i>x_j$; it is an \emph{antichain} if 
whenever $x_i\ne x_j$, then $x_i\not<x_j$ and $x_i\not>x_j$.
A chain $X$ is \emph{maximal} if whenever $y\in P\setminus X$, then $X\cup\{y\}$ is no longer a chain.
A poset $P$ is \emph{graded} if $P$ has the maximum and minimum elements and moreover, every maximal chain has the same length. The \emph{rank} of such $P$ is the length of a (any) maximal chain. Say two posets $(P, \le)$ and $(P', \le')$ are \emph{isomorphic} ($P\cong P'$) if
there is a bijection $f:P\to P'$ such that $x\le y$ $\iff$ $f(x)\le' f(y)$.

\subsection*{Lattice} \label{lat}

Let $P$ be as above (so that we deal only with \emph{finite} posets). Given $x, y\in P$, consider $\{z\in P \mid z\le x \mbox{ and } z\le y\}$.
If there exists a unique maximal element of this set, then we call it the \emph{meet} of $x$ and $y$ in $P$ (denoted by $x\wedge y$). We define the \emph{join} $x\vee y$ order-dually.
Say $P$ is a \emph{lattice} if $x\wedge y$ and $x\vee y$ exist for all $x, y\in P$.
 A subset $X\subseteq P$ is \emph{meet-dense} if whenever
$x\in P$, then there exists $Y\subseteq X$ such that $x=\wedge \,Y$.
Say $z$ is \emph{meet-irreducible} in $P$ if $z$ is not the maximum element of $P$ and whenever $z=x \wedge y$ then $z=x$ or $z=y$. Denote by $M(P)$ the set of all meet-irreducible elements in $P$.
A lattice $P$ is \emph{distributive} if 
$x\vee (y\we z)=(x\vee y)\we (x\vee z)$ and 
$x\we (y\vee z)=(x\we y)\vee (x\we z)$ for all $x, y, z$.
The \emph{Dedekind-MacNeille completion} of $P$ is the smallest lattice containing $P$.
Consequently, if $P\cong P'$ as posets, then their Dedekind-MacNeille completions are isomorphic.


\begin{figure}[t]
\caption{Dedekind-MacNeille completion} \label{ss3}
\[
\begin{xy}
0;<5mm,0mm>:
,(-3,2)*-{\ci}="al1"*++!R{}
,(3,2)*-{\ci}="ar1"*++!L{}
,(-3,5)*-{\ci}="cl1"*++!R{}
,(3,5)*-{\ci}="cr1"*++!L{}
,(5,3.5)*{\longrightarrow}
,(-3,2)+(10,0)*-{\ci}="al2"*++!R{}
,(0,3.5)+(10,0)*{\cci}="m2"*++!L{}
,(3,2)+(10,0)*-{\ci}="ar2"*++!L{}
,(-3,5)+(10,0)*-{\ci}="cl2"*++!R{}
,(3,5)+(10,0)*-{\ci}="cr2"*++!L{}
,(0,6.5)+(10,0)*{\cci}="t2"*++!L{}
,(0,0.5)+(10,0)*{\cci}="b2"*++!L{}
,\ar@{-}"al2";"b2"
,\ar@{-}"ar2";"b2"
,\ar@{-}"al2";"m2"
,\ar@{-}"ar2";"m2"
,\ar@{-}"cl2";"t2"
,\ar@{-}"cr2";"t2"
,\ar@{-}"cl2";"m2"
,\ar@{-}"cr2";"m2"
,\ar@{-}"al1";"cl1"
,\ar@{-}"ar1";"cr1"
,\ar@{-}"al1";"cr1"
,\ar@{-}"ar1";"cl1"
\end{xy}\]
\end{figure}
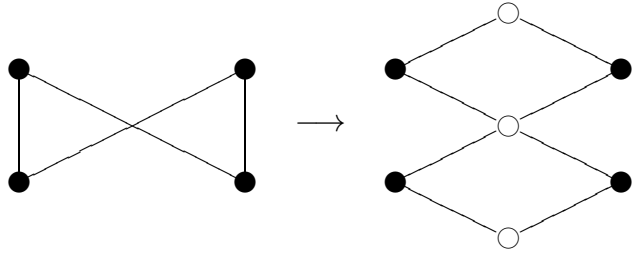


\bex{Figure \ref{ss3} illustrates an example of the Dedekind-MacNeille completion.
The poset on left is the original poset, and the one on right is its completion so that the meet and join all exist. Black dots $\ci$ indicate meet-irreducible elements.
}\eex

\bft{Every finite distributive lattice $P$ is a graded poset. 
Moreover, the function $m_P(x)=|\{z\in M(P) \mid x\le z\}|$
satisfies $x\lhd y \Longrightarrow$ $m_P(x)-m_P(y)=1$; this is the order dual of \cite[Proposition 2.10]{aigner}. In particular, the rank of $P$ is equal to $|M(P)|$.
}\eft


\subsection*{Coxeter group}\label{cox}

By $S_n$ we mean the symmetric group on $[n]:=\{1, 2, \dots, n\}$.
To represent permutations, we often use one-line notation: ``$w=i_1\cdots i_n$ with $i_k\in [n]$" means $w(k)=i_k$. For instance, $w=231$ means $w(1)=2, w(2)=3$ and $w(3)=1$.
The \emph{reverse permutation} is $i\mapsto n-i+1$.

For each $i\in [n-1]$, let $s_i$ denote the transposition interchanging $i$ and $i+1$. We call $S=\{s_i \mid i \in [n-1]\}$ \emph{Coxeter generators} of $S_n$. 
These elements satisfy \emph{Coxeter relations}:
\begin{align*}
s_i^2&=e,\\
s_{i}s_{i+1}s_i&=s_{i+1}s_{i}s_{i+1},\\
s_is_j&=s_js_i \mbox{  for $|i-j|\ge 2$.} 
\end{align*}
Moreover, $S$ is indeed a group-theoretic generator of $S_n$.
That is, for each $w\in S_n$, there exist $s_{i_1}, \dots, s_{i_k}$ such that
$w=s_{i_1}\cdots s_{i_k}$.
Let $\ell(w)=\min \{k\ge 0 \mid w=s_{i_1}\cdots s_{i_k}\}$ be the \emph{Coxeter length}.
Say $(i, j)$ is an inversion of $w$ if $i<j$ and $w(i)>w(j)$.
\bft{$\ell(w)$ is equal to the number of inversions of $w$.
}\eft


\subsection*{Bruhat order}\label{bru}
For $w\in S_n$ and $(i, j)\in [n]^2$, let $w(i, j)=|\{k  \mid k\le i \mbox{ and }w(k)\le j\}|$.
Define \emph{Bruhat order} $v\le w$ if $v(i, j)\ge w(i, j)$ for all $i, j$.
This gives a graded poset structure $(S_n, \le, \ell)$ of rank $n(n-1)/2$.
Also, define \emph{reverse Bruhat order} $v\le' w$ if $v(i, j)\le w(i, j)$ for all $i, j$.


\bft{$|M(S_n)|=(n-1)n(n+1)/6$.
}\eft
\begin{fact}
Let $x\in S_n$. Then 
$\beta(x):=\sum_{\substack{i<j\\x(i)>x(j)}}x(i)-x(j)$
is equal to the number of meet-irreducible permutations $z$ such that $x\le' z$ in reverse Bruhat order \cite[Theorem]{koba}.
\end{fact}

\bft{The Dedekind-MacNeille completion of $S_n$ is isomorphic to $\tn{ASM}_n$.
Moreover, $\tn{ASM}_n$ is a finite distributive lattice and $M(S_n)=M(\tn{ASM}_n)$.
}\eft

Hence we can unambiguously extend the function $\beta$ onto $\tn{ASM}_n$ which coincides with $m=m_{\tn{ASM}_n}$.

\bft{For all $x\in S_n$, we have
\[0\le \ell(x)\le \fr{n(n-1)}{2} \mbox{ and } 0\le \beta(x)\le \fr{(n-1)n(n+1)}{6}.\]
This is equivalent to 
\[-1\le1-\fr{2\ell(x)}{\fr{n(n-1)}{2}} \le 1 \mbox{ and } -1\le 1-\fr{2\beta(x)}{\fr{(n-1)n(n+1)}{6}}\le 1
\]
as Kendall's $\tau$ and Spearman's $\rho$ for sample distributions without ties.


}\eft


\bibliography{Kobayashi_copula}
\bibliographystyle{amsplain}

\end{document}